\documentclass[a4paper,11pt]{article}
\usepackage{amssymb,enumerate}
\usepackage{amsmath,amssymb}
\newtheorem{thm}{THEOREM}
\newtheorem{df}{DEFINITION}

\newtheorem{cor}[thm]{COROLLARY}
\newtheorem{rmk}{REMARK}
\def\bp{{\noindent\it Proof. \ }}

\numberwithin{equation}{section}

\title{\textbf{ATOMIC DECOMPOSITIONS FOR OPERATORS IN REPRODUCING KERNEL HILBERT SPACES}}

\author{{ LAURA G\u AVRU\c TA
\thanks{ The author was supported by the FWF project P 24986-N25}} \\
}
\date{}

\begin{document}

\maketitle

\begin{minipage}{120mm}
\small{\bf Abstract.} { In this paper, we give some results concerning atomic decompositions for operators on reproducing kernel Hilbert spaces, using frame theory techniques. We provide also generalizations for atomic decompositions of some theorems for reproducing kernel Hilbert spaces.\\
In particular, we obtain atomic decomposition results for operators on Bergman spaces and Fock spaces in a simple manner.}\\



\end{minipage}
\begin{center}
\section{\textbf{INTRODUCTION}}
\end{center}
In the following, we present some basic known facts about frames.
Frames were introduced by R.J. Duffin and A.C. Schaffer \cite{Duffin} in 1952, in the context of nonharmonic Fourier series. In 1986, frames were brought to life by Daubechies, Grossman and Meyer in the fundamental paper \cite{Daubechies2}.

We denote by $\mathcal{H}$ a separable Hilbert space and by $\mathcal{L}(\mathcal{H})$ the space of all linear bounded operators on $\mathcal{H}.$

 \begin{df}\label{df_frame} A family of elements $\{f_n\}_{n=1}^{\infty}\subset\mathcal{H}$ is called a frame of $\mathcal{H}$ if there exists
 constants $A,B>0$ such that
 $$A\|x\|^2\leq \sum_{n=1}^{\infty}|\langle x,f_n\rangle|^2\leq B\|x\|^2,\quad x\in\mathcal{H}.$$
 The constants $A,B$ are called frame bounds.
\end{df}

If just the last inequality in the above definition holds, we say that  $\{f_n\}_{n=1}^{\infty}$ is a \textit{Bessel sequence}.\\

The operator $$T:l^2\rightarrow \mathcal{H},\hspace{2mm}T\{c_n\}_{n=1}^{\infty}:=\sum_{n=1}^{\infty} c_nf_n$$ is called \textit{synthesis operator(or pre-frame operator)}. The adjoint operator is given by $$\Theta=T^{*}:\mathcal{H}\rightarrow l^2,\hspace{2mm}\Theta x=\{\langle x,f_n\rangle\}_{n=1}^{\infty}$$ and is called the \textit{analysis operator}. By composing $T$ with its adjoint $T^{*}$ we obtain the \textit{frame operator}
$$S:\mathcal{H}\rightarrow\mathcal{H},\hspace{2mm}Sx=TT^{*}x=\sum_{n=1}^{\infty}\langle x,f_n\rangle f_n.$$

The following Theorem is one of the most important result in frame theory because it gives the reconstruction formula.
\begin{thm} Let $\{f_n\}_{n=1}^{\infty}\subset \mathcal{H}$ be a frame for $\mathcal{H}$ with frame operator $S$. Then
\begin{enumerate}[($i$)]
\item $S$ is invertible and self-adjoint;
\item every $x\in \mathcal{H}$ can be represented as
\begin{equation}x=\sum_{n=1}^{\infty}\langle x,f_n\rangle  S^{-1}f_n=\sum_{n=1}^{\infty}\langle x,S^{-1}f_n\rangle f_n.
\end{equation}
\end{enumerate}
\end{thm}

For basic results on frame theory see the references \cite{Casazza}, \cite{Christensen}, \cite{Chui}, \cite{Daubechies}, \cite{Grochenig}, \cite{Heil}, \cite{Kovacevic}.\\

The results in this paper are organized as follows. In Section 2, we recall the definition of an atomic system for operators and we complete a Theorem for the characterization of these systems, considered for the first time by the author in \cite{LGavruta}. We rewrite this Theorem for a reproducing kernel Hilbert space. In Section 3 of this paper, we present some basic known facts about Bergman spaces, we give the atomic decomposition result for this case, for the standard weighted Bergman spaces and for Bergman spaces with B\'ekoll\'e-Bonami weights. In the last section of the paper, we give an atomic decomposition result for operators in Fock spaces.
\begin{center}
\section{\bf{ATOMIC DECOMPOSITIONS FOR OPERATORS}}
\end{center}
The first atomic decompositions results for holomorphic functions in Lebesgue spaces were obtained in 1980 by R.R. Coifman and R. Rochberg \cite{Coifman}.

Let $\mathcal{H}$ be a separable Hilbert space and $\{f_n\}_{n=1}^{\infty}\subset\mathcal{H}$.
\begin{df}\cite{LGavruta}\label{df_atomicsystem} Let $L\in\mathcal{L}(\mathcal{H}).$ We say that $\{f_n\}_{n=1}^{\infty}$ is an atomic system for $L$ if the following statements hold
\begin{enumerate}[($i$)]
\item the series $\displaystyle{\sum_n} c_nf_n$ converges for all $c=(c_n)\in l^2;$
\item there exists $C>0$ such that for every $x\in\mathcal{H}$ there exists $a_x=(a_n)\in l^2$ such that $\|a_x\|_{l^2}\leq C\|x\|$ and $Lx=\displaystyle{\sum_n} a_nf_n.$
\end{enumerate}
\end{df}
\begin{rmk}\label{Rmk0}
The condition $(i)$ in Definition \ref{df_atomicsystem} actually says that $\{f_n\}_{n=1}^{\infty}$ is a Bessel sequence. (see Corollary 3.2.4 in \cite{Christensen})
\end{rmk}

The above definition and the following main result (Theorem \ref{main_1}) presents a generalization of the concept of frames. The particular case of orthogonal projections was considered by  H.G. Feichtinger and T. Werther \cite{Feichtinger}.

The following Theorem for the existence of the atomic systems for an operator was proved in \cite{LGavruta}.
\begin{thm} Let $\mathcal{H}$ be a separable Hilbert space and $L\in\mathcal{L}(\mathcal{H})$. Then $L$ has an atomic system.
\end{thm}

In the following Theorem, we present a characterization for atomic systems.
\begin{thm}\label{main_1} Let $\{f_n\}_{n=1}^{\infty}\subset\mathcal{H}.$ Then the following statements are equivalent
\begin{enumerate}[($i$)]
\item $\{f_n\}_{n=1}^{\infty}$ is an atomic system for $L$;
\item there exists $A,B>0$ such that $$A\|L^{*}x\|^2\leq \sum_{n=1}^{\infty} |\langle x,f_n\rangle|^2\leq B\|x\|^2,\quad for \hspace{1mm} any\hspace{1mm}x\in\mathcal{H};$$
\item $\{f_n\}_{n=1}^{\infty}$ is a Bessel sequence and there exists a Bessel sequence $\{g_n\}_{n=1}^{\infty}$ such that $$Lx=\sum_{n=1}^{\infty} \langle x,g_n\rangle f_n;$$
\item $\displaystyle{\sum_{n=1}^{\infty}}|\langle x,f_n\rangle|^2<\infty$ and there exists a Bessel sequence $\{g_n\}_{n=1}^{\infty}$ such that $$L^*x=\sum_{n=1}^{\infty} \langle x,f_n\rangle g_n, \forall x\in\mathcal{H}.$$
\end{enumerate}
\end{thm}
\bp $(i)\Longleftrightarrow (ii)\Longleftrightarrow (iii)$ were proved in \cite{LGavruta}\\

$(iv)\Longrightarrow (ii)$ Using Cauchy-Schwartz inequality, we have:
\begin{align*}\langle L^*x, L^*x\rangle&=\langle \sum_{n=1}^{\infty}\langle x,f_n\rangle g_n, L^*x\rangle\\
&=\sum_{n=1}^{\infty}\langle x,f_n\rangle\langle g_n, L^*x\rangle\\
&=|\sum_{n=1}^{\infty}\langle x,f_n\rangle\langle g_n, L^*x\rangle|\\
&\leq\bigg(\sum_{n=1}^{\infty}|\langle x,f_n\rangle|^2\bigg)^{1/2}\bigg(\sum_{n=1}^{\infty}|\langle g_n,L^*x\rangle|^2\bigg)^{1/2}.
\end{align*}
Since $\{g_n\}_{n=1}^{\infty}$ is a Bessel sequence, we obtain
$$\langle L^*x, L^*x\rangle\leq C\|L^*x\|\bigg(\sum_{n=1}^{\infty}|\langle x,f_n\rangle|^2\bigg)^{1/2}.$$

If $L^*x\neq 0$, we have $\displaystyle \|L^*x\|^2\leq C\|L^*x\|\bigg(\sum_{n=1}^{\infty}|\langle x, f_n\rangle|^2\bigg)^{1/2}$ and from here we get $$\frac{1}{C^2}\|L^*x\|^2\leq\sum_{n=1}^{\infty}|\langle x, f_n\rangle|^2$$
From hypothesis, $\displaystyle \sum_{n=1}^{\infty}|\langle x, f_n\rangle|^2<\infty,\quad\forall x\in\mathcal{H}$ (actually $\{f_n\}_{n=1}^{\infty}$ is a Bessel sequence).
From \cite{Christensen} we have that there exists $B>0$ such that $$\displaystyle \sum_{n=1}^{\infty}|\langle x,f_n\rangle|^2\leq B\|x\|^2.$$\\

$(iii)\Longrightarrow (iv)$\\

From $(iii)$ we have $\displaystyle Lx=\sum_{n=1}^{\infty}\langle x,g_n\rangle f_n.$\\

 For any $y\in \mathcal{H}$ we have
 \begin{align*}
 \langle L^*y,x\rangle &=\langle y,Lx\rangle=\langle y,\sum_{n=1}^{\infty}\langle x,g_n\rangle f_n\rangle\\
                       &=\sum_{n=1}^{\infty}\langle g_n,x\rangle\langle y,f_n\rangle\\
                       &=\sum_{n=1}^{\infty}\langle \langle y,f_n\rangle g_n,x\rangle\\
                       &=\langle\sum_{n=1}^{\infty} \langle y,f_n\rangle g_n,x\rangle
 \end{align*}
 which implies $\displaystyle L^*y=\sum_{n=1}^{\infty}\langle y,f_n\rangle g_n\qquad\square$

\begin{df}\label{dfLframe}
We say that $\{f_n\}_{n=1}^{\infty}$ is a $L-$frame (or a frame for $L$) if there exists the constants $A,B>0$ such that $$A\|L^*x\|^2\leq\sum_{n=1}^{\infty}|\langle x,f_n\rangle|^2\leq B\|x\|^2,\quad\forall x\in\mathcal{H}.$$
\end{df}

We consider now the reproducing kernel Hilbert space $(\mathcal{H},K)$, with the kernel $K(z,\lambda)=K_{\lambda}(z)$, $K:\Omega\times\Omega\rightarrow\mathbb{C}$, where $\Omega$ is a given nonempty set.

We recall that a reproducing kernel Hilbert space $(\mathcal{H},K)$ is a Hilbert space $\mathcal{H}$ of functions on $\Omega$ such that for every $\lambda$, $K_{\lambda}$ belongs to $\mathcal{H}$ and for every $\lambda\in\Omega$ and every $f\in\mathcal{H},$ $f(\lambda)=\langle f,K_{\lambda}\rangle$.

Let $\Lambda=\{\lambda_n\}_{n=1}^{\infty}\subset\Omega$ be a subset of points with $\lambda_n\neq\lambda_m$, for $n\neq m$.
In this case, with $f_n=\dfrac{K_{\lambda_n}}{\|K_{\lambda_n}\|}$, $n=1,2,\ldots$, we rewrite Theorem \ref{main_1} as follows:
\begin{thm}\label{rkhs-g} The following statements are equivalent
\begin{enumerate}[($i$)]
\item $\bigg\{\dfrac{K_{\lambda_n}}{\|K_{\lambda_n}\|}\bigg\}_{n=1}^{\infty}$ is an atomic system for $L$ i.e.\\
\begin{enumerate}[$\circ$]
\item the series $\displaystyle \sum_{n=1}^{\infty} c_n\dfrac{K_{\lambda_n}}{\|K_{\lambda_n}\|}$ converges for all $\{c_n\}\in l^2$; and
    \item there exists $C>0$ such that for every $f\in\mathcal{H}$ there exists $a_f=(a_n)\in l^2$ such that $\|a_f\|_{l^2}\leq C\|f\|$ and $\displaystyle Lf=\sum_{n=1}^{\infty} a_n\dfrac{K_{\lambda_n}}{\|K_{\lambda_n}\|}.$
\end{enumerate}
\item there exists $A,B>0$ such that $$A\|L^{*}f\|^2\leq \sum_{n=1}^{\infty} \frac{|f(\lambda_n)|^2}{\|K_{\lambda_n}\|^2}\leq B\|f\|^2,\quad for \hspace{1mm} any\hspace{1mm}f\in\mathcal{H};$$
\item $\bigg\{\dfrac{K_{\lambda_n}}{\|K_{\lambda_n}\|}\bigg\}_{n=1}^{\infty}$ is a Bessel sequence and there exists a Bessel sequence $\{g_n\}_{n=1}^{\infty}$ such that $$Lf=\sum_{n=1}^{\infty} \langle f,g_n\rangle \frac{K_{\lambda_n}}{\|K_{\lambda_n}\|};$$
\item $\displaystyle{\sum_{n=1}^{\infty}}\frac{|f(\lambda_n)|^2 }{\|K_{\lambda_n}\|^2}<\infty$ and there exists a Bessel sequence $\{g_n\}_{n=1}^{\infty}$ such that $$L^*f=\sum_{n=1}^{\infty} \dfrac{f(\lambda_n)}{\|K_{\lambda_n}\|}g_n, \forall f\in\mathcal{H}.$$
\end{enumerate}
\end{thm}

\begin{df}We say that a sequence of distinct points $\{\lambda_n\}_{n=1}^{\infty}$ in $\Omega$ is a sampling sequence for $L$ if there exists two positive constants $A, B$ such that \begin{equation}\label{samplingsequence}A\|L^*f\|^2_{\mathcal{H}}\leq\sum_{n=1}^{\infty}\frac{| f(\lambda_n)|^2}{\|K_{\lambda_n}\|^2}\leq B\|f\|^2_{\mathcal{H}}\end{equation}
\end{df}
In the case when $L$ is the identity operator, we have the notion of sampling sequence for $\mathcal{H}.$
The normalized kernel is given by the following relation: $$k_{\lambda}(z):=\frac{K_{\lambda}(z)}{\|K_{\lambda}\|}.$$
\begin{rmk}The relation (\ref{samplingsequence}) in Theorem \ref{rkhs-g} is equivalent with the fact that $\{k_{\lambda_n}\}_{n=1}^{\infty}$ is a frame for $L$: $$A\|L^*f\|^2_{\mathcal{H}}\leq\sum_{n=1}^{\infty}|\langle f,k_{\lambda_n}\rangle|^2\leq B\|f\|^2_{\mathcal{H}}$$ because
\begin{align*}\langle f,k_{\lambda_n}\rangle &=\langle f,\dfrac{K_{\lambda_n}}{\|K_{\lambda_n}\|_{\mathcal{H}}}\rangle\\
&=\dfrac{1}{\|K_{\lambda_n}\|_{\mathcal{H}}}\langle f, K_{\lambda_n}\rangle\\
&=\dfrac{1}{\|K_{\lambda_n}\|_{\mathcal{H}}}f(\lambda_n).
\end{align*}
\end{rmk}
\begin{center}
\section{\textbf{ATOMIC DECOMPOSITIONS FOR OPERATORS IN BERGMAN SPACES}}
\end{center}
Let $\mathbb{C}$ be the complex plane and let $\mathbb{D}=\{z\in\mathbb{C}:|z|<1\}$ be the open unit disk in $\mathbb{C}.$ We denote by $dA(z)$ the area measure on $\mathbb{D}$, normalized such that the area measure on $\mathbb{D}$ is $1:$
$$dA(z)=\frac{1}{\pi}dxdy=\frac{1}{\pi}rdrd\theta,\quad z=x+iy=re^{i\theta}.$$

The Bergman space $L_a^2$ (sometimes denoted by $A^2$) is the subset of $L^2(\mathbb{D})$ consisting of analytic functions. $L_a^2$ is a reproducing kernel Hilbert space, with reproducing kernel : $$K(z,\lambda)=K_\lambda(z)=\frac{1}{(1-\overline{\lambda}z)^2}.$$
and normalized kernel $$k_{\lambda}(z)=\dfrac{K_{\lambda}(z)}{\|K_{\lambda}\|}=\dfrac{1-|\lambda|^2}{(1-\overline{\lambda}z)^2},\quad \lambda\in\mathbb{D}.$$

For a detailed account results of atomic decompositions for the identity operator on Bergman spaces see \cite{Zhu}.

Then, from Theorem \ref{rkhs-g}, we obtain the following Corollary:
\begin{cor}Let $\{\lambda_n\}_{n=1}^{\infty}\subset \mathbb{D}$. The following are equivalent
\begin{enumerate}[($i$)]
\item $\{k_{\lambda_n}\}_{n=1}^{\infty}$ is an atomic system for $L$ i.e.
\begin{enumerate}[$\circ$]
\item the series $\displaystyle \sum_{n=1}^{\infty} c_n\dfrac{1-|\lambda_n|^2}{(1-\overline{\lambda}_n z)^2}$ converges in $L_a^2$ for all $\{c_n\}\in l^2$; and
    \item there exists a positive constant $C$ such that for every\\ $f\in\mathcal{H}$ there exists $a_f=(a_n)\in l^2$ such that $\|a_f\|_{l^2}\leq C\|f\|$ and $\displaystyle Lf(z)=\sum_{n=1}^{\infty} a_n\dfrac{1-|\lambda_n|^2}{(1-\overline{\lambda}_n z)^2};$
\end{enumerate}
\item there exists $A, B>0$ such that $\{\lambda_n\}_{n=1}^{\infty}$ is a sampling sequence for $L;$ i.e. there exists two positive constants $A, B$ such that \begin{equation*}\label{samplingsequence}A\|L^*f\|^2_{\mathcal{H}}\leq\sum_{n=1}^{\infty}|f(\lambda_n)|^2(1-|\lambda_n|^2)^2\leq B\|f\|^2_{\mathcal{H}};\end{equation*}
\item $\{k_{\lambda_n}\}_{n=1}^{\infty}$ is a Bessel sequence and there exists a Bessel sequence $\{g_n\}_{n=1}^{\infty}$ such that $$Lf=\sum_{n=1}^{\infty}\langle f,g_n\rangle \dfrac{1-|\lambda_n|^2}{(1-\overline{\lambda}_n z)^2};$$
    \item $\displaystyle \sum_{n=1}^{\infty}(1-|\lambda_n|^2)^2|f(\lambda_n)|^2<\infty$ and there exists a Bessel sequence $\{g_n\}_{n=1}^{\infty}$ such that $$L^*f=\sum_{n=1}^{\infty}(1-|\lambda_n|^2)f(\lambda_n)g_n.$$
\end{enumerate}
\end{cor}

More general, for $\eta>-1$ let $dA_{\eta}$ be the area measure on $\mathbb{D}$, $$dA_{\eta}(z)=(\eta+1)(1-|z|^2)^{\eta}dA(z).$$

The standard weighted Bergman space $L_a^2(dA_{\eta})$ is the subset of $L^2(\mathbb{D},dA_{\eta})$, consisting of analytic functions.
  $$\displaystyle L^2(\mathbb{D},dA_{\eta})=\{f:\mathbb{D}\rightarrow\mathbb{C}: \|f\|_{2,\eta}^2:=\int_{\mathbb{D}}|f(z)|^2dA_{\eta}<\infty\}.$$
 Then $L_a^2(dA_{\eta})$ is a RKHS with reproducing kernel $$K^{\eta}(z,\lambda)=K_{\lambda}^{\eta}(z)=\frac{1}{(1-\overline{\lambda}z)^{2+\eta}},\quad z,\lambda\in\mathbb{D}$$
and normalized reproducing kernel $$k_{\lambda}^{\eta}(z)=\frac{K_{\lambda}^{\eta}(z)}{\|K_{\lambda}^{\eta}\|_{2,\eta}}=\frac{(1-|\lambda|^2)^{1+\frac{\eta}{2}}}{(1-\overline{\lambda}z)^{2+\eta}},\quad z,\lambda\in\mathbb{D}.$$

\begin{cor}Let $\{\lambda_n\}_{n=1}^{\infty}\subset \mathbb{D}$. The following are equivalent
\begin{enumerate}[($i$)]
\item $\{k_{\lambda_n}^{\eta}\}_{n=1}^{\infty}$ is an atomic system for $L$ i.e.
\begin{enumerate}[$\circ$]
\item the series $\displaystyle \sum_{n=1}^{\infty} c_n\dfrac{(1-|\lambda_n|^2)^{1+\frac{\eta}{2}}}{(1-\overline{\lambda}_nz)^{2+\eta}}$ converges in $L_a^2(dA_{\eta})$ for all\\ $\{c_n\}\in l^2$; and
    \item there exists $C>0$ such that for every $f\in\mathcal{H}$ there exists\\ $a_f=(a_n)\in l^2$ such that $\|a_f\|_{l^2}\leq C\|f\|$ and $\displaystyle Lf(z)=\sum_{n=1}^{\infty} a_n\dfrac{(1-|\lambda_n|^2)^{1+\frac{\eta}{2}}}{(1-\overline{\lambda}_nz)^{2+\eta}};$
\end{enumerate}
\item there exists $A, B>0$ such that $\{\lambda_n\}_{n=1}^{\infty}$ is a sampling sequence for $L;$ i.e. there exists two positive constants $A,B$ such that $$A\|L^*f\|^2\leq\sum_{n=1}^{\infty}|f
    (\lambda_n)|^2(1-|\lambda_n|^2)^{2+\eta}\leq B\|f\|^2$$
\item $\{k_{\lambda_n}^{\eta}\}_{n=1}^{\infty}$ is a Bessel sequence and there exists a Bessel sequence $\{g_n\}_{n=1}^{\infty}$ such that $$Lf=\sum_{n=1}^{\infty}\langle f,g_n\rangle \dfrac{(1-|\lambda_n|^2)^{1+\frac{\eta}{2}}}{(1-\overline{\lambda}_nz)^{2+\eta}};$$
    \item $\displaystyle \sum_{n=1}^{\infty}(1-|\lambda_n|^2)^{2+\eta}|f(\lambda_n)|^2<\infty$ and there exists a Bessel sequence $\{g_n\}_{n=1}^{\infty}$ such that $$L^*f=\sum_{n=1}^{\infty}(1-|\lambda_n|^2)^{1+\frac{\eta}{2}}f(\lambda_n)g_n.$$
\end{enumerate}
\end{cor}
\begin{rmk} The above Corollary is a generalization of some results in \cite{Vaezpour}.
\end{rmk}

 The weighted Bergman space $A^2(\omega)=L_a^2(\omega)$ is the subset $L^2(\omega)$, consisting of analytic functions.
  $$L^2(\omega)=\{f:\mathbb{D}\rightarrow\mathbb{C}:\|f\|_{A^2(\omega)}^2:=\int_{\mathbb{D}}|f(z)|^2\omega(z)dA(z)<\infty\}.$$

  For $\eta>-1$ the class $B_2(\eta)$ consists of weights $\omega$  with the property that there exists a constant $c>0$ such that
$$\bigg(\int_{S(\theta,h)}\omega dA_{\eta}\bigg)\bigg(\int_{S(\theta,h)}\omega^{-1} dA_{\eta}\bigg)\leq c[A_{\eta}(S(\theta,h))]^2$$
for any Carleson square: $$S(\theta,h)=\bigg\{z=re^{i\alpha}:1-h<r<1, |\theta-\alpha|<\dfrac{h}{2}\bigg\},\quad \theta\in[0,2\pi], h\in(0,1).$$
The weights $\omega$ considered here are called B\'ekoll\'e weights.\\

In the following, we consider $\alpha\in(0,1)$, $\eta>-1$ and $\dfrac{\omega}{(1-|z|^2)^\eta}\in B_2({\eta}).$ Then $L_a^2(\omega)$ is a reproducing kernel Hilbert space. We denote by $$K^{\eta,\omega}(z,\lambda)=K_{\lambda}^{\eta,\omega}(z)$$ the reproducing kernel of this space.\\

 In \cite{OConstantin}, O. Constantin gave the following estimation for the norm of Bergman kernel $K_{\lambda}^{\eta,\omega}:$
 $$\|K_{\lambda}^{\eta,\omega}\|^2\sim\bigg(\int_{D_{\lambda,\alpha}}\omega dA\bigg)^{-1},$$ for every disc $D_{\lambda,\alpha}=\{z\in\mathbb{D}:|z-\lambda|<\alpha(1-|\lambda|)\}.$\\ By $\sim$ we mean that the involved constants are independent of $\lambda\in\mathbb{D}.$ (For two real valued functions $E_1$, $E_2$ we write $E_1\sim E_2$ if there there exists a positive constant $c$ independent of the argument such that $\dfrac{1}{c}E_1\leq E_2\leq c E_1$).

\begin{cor}
Let $\{\lambda_n\}_{n=1}^{\infty}\subset\mathbb{D}.$ Then the following statements are equivalent
\begin{enumerate}[($i$)]
\item $\displaystyle\bigg\{K_{\lambda_n}^{\eta,\omega}\bigg(\int_{D_{\lambda_n,\alpha}}\omega dA\bigg)^{1/2}\bigg\}_{n=1}^{\infty}$ is an atomic system for $L$, i.e.
\begin{enumerate}[$\circ$]
\item the series $\displaystyle \sum_{n=1}^{\infty} c_n\bigg\{K_{\lambda_n}^{\eta,\omega}\bigg(\int_{D_{\lambda_n,\alpha}}\omega dA\bigg)^{1/2}\bigg\}$
 converges in $L_a^2(\omega)$ for all $\{c_n\}\in l^2$; and
    \item there exists $C>0$ such that for every $f\in\mathcal{H}$ there exists\\ $a_f=(a_n)\in l^2$ such that $\|a_f\|_{l^2}\leq C\|f\|$\\ and $\displaystyle Lf(z)=\sum_{n=1}^{\infty} a_n\bigg\{K_{\lambda_n}^{\eta,\omega}\bigg(\int_{D_{\lambda_n,\alpha}}\omega dA\bigg)^{1/2}\bigg\};$
\end{enumerate}
\item there exists $A, B>0$ such that $\{\lambda_n\}_{n=1}^{\infty}$ is a sampling sequence for $L;$ i.e. there exists two positive constants $A, B$ such that $$A\|L^*f\|^2_{\mathcal{H}}\leq\sum_{n=1}^{\infty}|f(\lambda_n)|^2\bigg(\int_{D_{\lambda_n,\alpha}}\omega dA\bigg)\leq B\|f\|^2_{\mathcal{H}};$$
\item $\displaystyle\bigg\{K_{\lambda_n}^{\eta,\omega}\bigg(\int_{D_{\lambda_n,\alpha}}\omega dA\bigg)^{1/2}\bigg\}_{n=1}^{\infty}$ is a Bessel sequence and there exists a Bessel sequence $\{g_n\}_{n=1}^{\infty}$ such that $$Lf=\sum_{n=1}^{\infty}\langle f,g_n\rangle\bigg\{K_{\lambda_n}^{\eta,\omega}\bigg(\int_{D_{\lambda_n,\alpha}}\omega dA\bigg)^{1/2}\bigg\}; $$
    \item $\displaystyle\sum_{n=1}^{\infty}|f(\lambda_n)|^2\bigg(\int_{D_{\lambda_n,\alpha}}\omega dA\bigg)<\infty$ and there exists a Bessel sequence $\{g_n\}_{n=1}^{\infty}$ such that $$L^*f=\sum_{n=1}^{\infty}f(\lambda_n)\bigg(\int_{D_{\lambda_n,\alpha}}\omega dA\bigg)^{1/2}g_n.$$
\end{enumerate}
\end{cor}

The above Corollary provides a generalization of Theorem 2.3 in \cite{Chacon}. Also, our proof is more simple, with a different technique.
\begin{center}
\section{ATOMIC DECOMPOSITIONS FOR OPERATORS IN FOCK SPACES}
\end{center}
Let $dA(z)$ be the usual Lebesgue measure on $\mathbb{C}$. For $\alpha>0$, $dA_{\alpha}(z)$ is defined as follows $$dA_{\alpha}(z)=\frac{\alpha}{\pi}e^{-\alpha|z|^2}dA(z).$$
The Fock space $F_{\alpha}^2$ is the space of all entire functions $f$ on $\mathbb{C}$ for which
$$\|f\|_{F_\alpha^2}^2=\int_{\mathbb{C}}|f(z)|^2dA_\alpha(z)<\infty.$$\\
 $F_\alpha^2$ is a RKHS, with reproducing kernel $$K(z,\lambda)=K_{\lambda}(z)=e^{\alpha z\overline{\lambda}}$$
and normalized reproducing kernel $$k_{\lambda}^{\alpha}(z)=\frac{K_{\lambda}(z)}{\|K_{\lambda}\|}=\frac{e^{\alpha z|\overline{\lambda}}}{e^{\frac{\alpha |\lambda|^2}{2}}}=
e^{\alpha(z\overline{\lambda}-\frac{|\lambda|^2}{2})}.$$

For details on Fock spaces see the book of K. Zhu \cite{Zhu1}. See also \cite{Tung}.
Atomic decompositions for the identity operator on Fock spaces were obtained in \cite{Janson}.
\begin{cor}Let $\{\lambda_n\}_{n=1}^{\infty}\subset \mathbb{D}$. The following are equivalent
\begin{enumerate}[($i$)]
\item $\{k_{\lambda_n}^{\alpha}\}_{n=1}^{\infty}$ is an atomic system for $L$ i.e.
\begin{enumerate}[$\circ$]
\item the series $\displaystyle \sum_{n=1}^{\infty} c_n e^{\alpha(z\overline{\lambda}_n-\frac{|\lambda_n|^2}{2})}$ converges in $F_{\alpha}^2$ for all $\{c_n\}\in l^2$; and
    \item there exists $C>0$ such that for every $f\in\mathcal{H}$ there exists\\ $a_f=(a_n)\in l^2$ such that $\|a_f\|_{l^2}\leq C\|f\|$ and $\displaystyle Lf(z)=\sum_{n=1}^{\infty} a_n e^{\alpha(z\overline{\lambda}_n-\frac{|\lambda_n|^2}{2})} .$
\end{enumerate}
\item there exists $A, B>0$ such that $\{\lambda_n\}_{n=1}^{\infty}$ is a sampling sequence for $L;$ i.e. there exists two positive constants $A,B$ such that $$A\|L^*f\|^2\leq\sum_{n=1}^{\infty}|f(\lambda_n)|^2e^{-\alpha|\lambda_n|^2}\leq B\|f\|^2;$$
\item $\{k_{\lambda_n}^{\alpha}\}_{n=1}^{\infty}$ is a Bessel sequence and there exists a Bessel sequence $\{g_n\}_{n=1}^{\infty}$ such that $$Lf=\sum_{n=1}^{\infty}\langle f,g_n\rangle  e^{\alpha(z\overline{\lambda}_n-\frac{|\lambda_n|^2}{2})};$$
    \item $\displaystyle \sum_{n=1}^{\infty}e^{-\alpha|\lambda_n|^2}|f(\lambda_n)|^2<\infty$
    and there exists a Bessel sequence $\{g_n\}_{n=1}^{\infty}$ such that $$L^*f=\sum_{n=1}^{\infty}e^{-\alpha\frac{|\lambda_n|^2}{2}}f(\lambda_n)g_n.$$
\end{enumerate}
\end{cor}
\begin{rmk} The above Corollary is a generalization of some results in \cite{Tatari}.
\end{rmk}
\textbf{Acknowledgements.} I would like to thanks Dr. Olivia Constantin for
introducing me in the theory of  Bergman and Fock spaces and also for the careful reading of this paper and her useful comments.
\newpage
\begin{center}

\end{center}
\begin{flushright}
\textit{
{ \normalsize Politehnica University of Timi\c soara,\\ Department of Mathematics }\\
{ \normalsize Pia\c ta Victoriei no. 2, 300006 Timi\c soara, Romania}\\
}
\vspace{3mm}
\textit{
{ \normalsize Fakult\"at fur Mathematik, University of Viena}\\
{ \normalsize Oskar-Morgenstern-Platz 1
1090 Wien, Austria}}\\

\hspace{5mm}{ \normalsize \textit{E-mail}: gavruta\_laura@yahoo.com}
\end{flushright}

\end{document}